\documentclass{ws-ijmpe}
\usepackage[ansinew]{inputenc}
\usepackage{fancyvrb}   
\usepackage{float}
\usepackage{graphicx}

\floatstyle{ruled}  
\newfloat{codigo}{tbp!}{lop}
\floatname{codigo}{Code}

\begin{document}

\markboth{B. O. Rodrigues, L. A. C. P. da Mota and L. G. S. Duarte}{Numerical Calculation With Arbitrary Precision}

\catchline{}{}{}{}{}

\title{NUMERICAL CALCULATION WITH ARBITRARY PRECISION}

\author{\footnotesize BRUNO OSÓRIO RODRIGUES, L. A. C. P. DA MOTA and L. G. S. DUARTE}
\address{Departamento de Física Teórica, Universidade do Estado do Rio de Janeiro, Rua S\~ao Francisco Xavier 524\\
Rio de janeiro, RJ 20550-900, Brazil\\
brunooz@gmail.com}

\maketitle

\begin{history}
\received{(14 May 2007)}
\revised{(21 Jun 2007)}
\accepted{(21 Jun 2007)}
\end{history}

\begin{abstract}
 The vast use of computers on scientific numerical computation makes  the
 awareness of the limited precision that these machines are able to
 provide us an essential matter.  A limited and insufficient precision
 allied to the truncation and rounding
 errors may induce the user to  incorrect interpretation of his/hers answer.
 In this work, we have developed a computational package to minimize this
 kind of error by offering arbitrary precision numbers and calculation.
 This is very important in Physics where we can work with numbers too
 small and too big simultaneously.
 \end{abstract}

\section{Introduction}
\begin{figure}[th]
\vspace*{30pt}
\includegraphics[width=1.5cm,height=1.2cm]{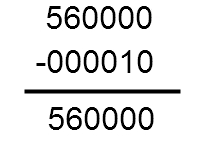}
\vspace*{12pt} \caption{An example of a wrong answer if were used
numbers with four digits of precision.} \label{fig1}
\end{figure}

\hspace\parindent In the past, when doing calculations by hand, on paper, we imagined that the numbers we were using had the precision we wanted, i.e., virtually unlimited. When we wrote ``2'', were actually regarding  it as meaning ``$2.0 \cdots 0$'',  with and infinity of zeros.  Unfortunately, life is not that simple and the numbers in the computer can not have as much digits as desired. It is only able to represent numbers with a limited length. Today, it is usual to work with numbers which present 32bits precision, this allows for calculations with sixteen digit numbers. So, one can conclude that the computer is not able to represent real numbers. This limited precision of the computers may induce scientists to wrong answers, and worse: they may be using those wrong answers as if they were the right ones (see fig. \ref{fig1}) .

There are some ways to handle this problems, the obvious one is to use numbers as large as needed.  In this undergraduate project, we are developing a computational package that permits us to make calculations with arbitrary precision numbers.


\section{A Solution: NOz}
\hspace\parindent NOz is the name of the computational package being developed\footnote{Besides being developed, there is already a functional version of NOz} in UERJ. It defines a new float point numerical type of arbitrary precision and the basic operations with them. Basically, it means that the user has the freedom to choose the precision of the numbers higher than the usual 32bits (or even smaller than it).

We know that there are already computational packages that can handle this problem, like GMP, but they are not easy to work with and demand some training before its use.\footnote{GNU Multi-Precision - An open source package, written in C++ and assembler. More details in its web site http://gmplib.org/ } Actually, the main feature that NOz has, as we shall explain later, is that it is very friendly.  Just declare it in your main program and use it.

This package is written in pascal and in further versions, it will be compatible with C++ and Fortran programs.

\section{The NOz Numerical Type}
\begin{figure}[th]
\vspace*{10pt}
\includegraphics[width=1.3cm,height=1.2cm]{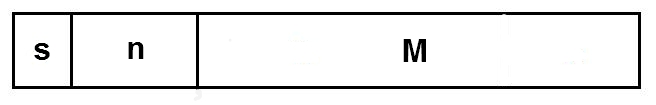}
\caption{The structure of a real number.} \label{fig2}
\end{figure}

\hspace\parindent The NOz arbitrary precision type is a structure (see fig. \ref{fig2}) that contains the necessary informations of a float point number: signal (s), mantissa(M) and exponent(n).

In NOz, the signal is represented by a boolean type, the exponent is an integer type and the mantissa is a dynamic array of natural numbers. Each position of the mantissa array represents a digit of the number. It means that the mantissa array will have the same number of positions as the number of digits in the number. The positions are filled from the less significant digit to the most significant one. The complexity of the NOz algorithms is proportional to the number of positions in the mantissa. We are currently developing better algorithms in order to be able to use more than one digit per position in the mantissa array to improve the speed of the calculations.

The precision, in NOz, is not a property of each number in a calculation. It is defined as a variable in the main unit of the package and passed on as a parameter to the arbitrary precision functions. It is a way to ensure that, in an arithmetic operation, the numbers involved have the same precision. As it is a variable, the user has the freedom to change the precision of the calculations in runtime.

\section{Numerical Operations Already Supported}
\hspace\parindent It's already possible to use the following operations involving numbers of arbitrary precision using NOz: sum, subtraction, multiplication, division, comparison operations (equal to, minor than, greater than, etc) and the factorial.

With the four basic arithmetic operations, all the polynomial problems can be dealt with. Functions such as exponential, sin, cos and others can be treated by using series as the factorial function is already implemented.

\section{Main Characteristics of NOz and Usage}
\hspace\parindent The NOz (will) have the follow characteristics:
\begin{itemlist}
 \item Supports numbers with precision up to 2 billion digits.
 \item Easy to use.
 \item Allows the use of the operator symbols ($+$, $-$, $*$ and $/$) by operator overload in the languages where it is possible, for instance, like Pascal.
\end{itemlist}

To use NOz, the user has just to declare the unoz.pas unit in its main program, set the value of precision variable in unoz unit and ensure that the archive noz.dll (the compiled NOz package) is in the same directory of the main program.
Code 1 is a little example of a program that calculates the exponential $e^x$ using its Taylor series. The functions \textit{sum}, \textit{divis}, \textit{multi} and \textit{factorial} are the arbitrary precision functions for sum, division, multiplication and factorial present in NOz.

Using $x = 1$, $n = 70$ and the precision set to 100 digits, the program results: $e = $ 2.7182818284590452353602874713526624977572470936999595749669676277 24076630353547594571382178525166427E0.

\begin{codigo}[h]
      \caption{Exponential}
      \vspace{5mm}
      \small
\begin{Verbatim}
program Exponential_NOz;
{$APPTYPE CONSOLE}
uses
   FastShareMem, SysUtils, UNOz {the NOz interface Unit};

function exponential (x : NOz_Float; n:integer) : NOz_Float;
var i : integer;
    e_x : NOz_Float;
begin
   e_x := NOz_0;

   for i := 0 to n do
      e_x := sum(e_x,divis(multi(x,x),factorial(i)));

   result := e_x;
end;

var x : NOz_Float;  // Declares x as an arbitrary precision float number
    n : integer;
begin
   x := StrToNoz('1.0E0'); //Convertion Routine from String to NOz_Float
   n := 70;  //Number of steps in the Exponential series

   writeln(NozToStr(exponential(x,n))); //Converts result to String
end.
\end{Verbatim}
   \label{code:expo}
\end{codigo}

\section{Conclusion}
\hspace\parindent Being conscious of the limited precision problem is as important as being able to use arbitrary precision numbers. The NOz is a good  solution to minimize the problem. Again, the problem is always there, we have to control it through using the NOz numbers that, although not being real numbers, present arbitrary precision, i.e., we can use the precision needed for the job at hand. The present version of NOz can deal with polynomial calculations only. But this is already of great use to the community since many interesting problems are polynomial in essence. Further versions will bring elementary functions as built-in operations. The main reason of its development is to give the community an easier way to handle the precision problem. But the use of  arbitrary precision routines is indicated only when it is really necessary since, as we increase the precision of the calculation, we have a corresponding decreasing of the speed of the calculation.

\section*{Acknowledgements}

Thanks to CNPq, FAPERJ and UERJ for the financial support and special thanks to Anibal L. Pereira, Luiz Fernando de Oliveira, J. Avellar and M. E. Bracco.


\begin{thebibliography}{0}

  \bibitem{digitais} G.G. Langdon Jr, E. Fregni, \textit{Projeto de Computadores Digitais} (Edgard Blücher, São Paulo, 1977)

  \bibitem{cschaum} John R. Hubbard, \textit{Programming With C++} (McGraw-Hill, New York, 1996)

  \bibitem{computerarit} Behrooz Parhami, \textit{Computer Arithmetic} (Oxford University Press, New York, 2000)

    \bibitem{bibliadelphi} Marcus Cantù, {\it Dominando o Delphi 6 ,``A Bíblia''} (Makron Books, São Paulo, 2002)

\end{thebibliography}
\end{document}